\documentclass[american]{article}
\usepackage{babel,amssymb,amsfonts,amscd,xypic}
\usepackage{amsmath}
\usepackage{times}
\newtheorem{theorem}{Theorem}
\newtheorem{proposition}{Proposition}

\newtheorem{defn}{Definition}

\everymath = {\rm}

\begin{document}

\title{Noncommmutative theorems: Gelfand Duality, Spectral, Invariant Subspace, and Pontryagin Duality}
\author{Mukul S. Patel}
\date{February 18, 2005}
\maketitle

\begin{abstract}
We extend the Gelfand-Naimark duality of commutative $C^*$-algebras,
\begin{center}
\scriptsize{\{A COMMUTATIVE $C^*$-ALGEBRA\} ------ \{A LOCALLY COMPACT HAUSDORFF SPACE\}}\newline
to\newline
\{A $C^*$-ALGEBRA\} ------ \{A QUOTIENT OF A LOCALLY COMPACT HAUSDORFF SPACE\}.
\end{center}
Thus, a $C^*$-algebra is isomorphic to the convolution algebra of continuous regular Borel measures on the topological equivalence relation given by the above mentioned quotient. In commutative case  this reduces to Gelfand-Naimark theorem. Applications: 1) A simultaneous extension, to arbitrary Hilbert space operators, of Jordan Canonical Form and  Spectral Theorem of normal operators  2) A functional calculus for arbitrary operators.  3) Affirmative solution of Invariant Subspace Problem. 4) Extension of Pontryagin duality to nonabelian groups, and inevitably to groups whose underlying topological space is noncommutative.
\end{abstract}

\setcounter{secnumdepth}{4}
\setcounter{tocdepth}{1}

\normalsize

\section{\normalsize{INTRODUCTION}}

Connes attaches $C^*$-algebras to various quotient spaces arising in geometry \cite{connes}. Conversely, we assign a natural quotient space to any given $C^*$-algebra. Of course, for commutative algebras,
the Gelfand-Naimark theorem does the job:
\begin{theorem}[Gelfand-Naimark]\label{gelfand}
A commutative $C^*$-algebra  $A$ is naturally isomorphic to $C_0(P(A))$, the algebra of complex valued continuous functions vanishing at infinity, on the (locally compact Hausdorff) maximal ideal space $P(A)$ of $A$.
\end{theorem}
Since the original theorem \cite{gelfand}, there have been several noncommutative generalizations in various directions \cite{1,2,5,4,7,enock} with varying degree of success. Our generalization (Theorem \ref{main})  is implemented by identifying the \emph{natural} noncommutative analog of locally compact Hausdorff space---a quotient of a locally compact Hausdorff space---and as such,  embraces most commonly ocurring geometric situations on one hand, and all $C^*$-algebras on the other. The key to this quotient is the following trivial observation: \textit{A $C^*$-algebra is commutative if and only if
all its irreducible Gelfand-Naimark-Segal representations
are pair-wise inequivalent.} Thus, the noncommutativity of an algebra is completely captured
by the equivalence relation given by equivalence of irreducible GNS representations.

There have been studies of $C^*$-algebras via continuous functions on groupoids (See \cite{renault}, for example). The latter include equivalence relations as a special case. However, a larger algebra is needed to capture the whole situation. Our main result, Theorem \ref{main}, asserts that the algebra $A$ is canonically isomorphic to a certain algebra of regular Borel \emph{measures} on an equivalence relation $R(A)$. We then take this equivalence relation, or equivalently the quotient map it entails, as quantum space.

\begin{defn}[Quantum spaces]\label{quantum space} 
A \textbf{quantum space (resp. compact quantum space)} is a quotient map $q: X\twoheadrightarrow Y$ where $X$ is a locally compact (resp. compact) Hausdorff space.  A \textbf{quantum group (resp. semigroup, groupoid, etc.)} is a group object (resp. semigroup object, groupoid object, etc.) in the category of quantum spaces. In this setting, the terms \textnormal{`abelian'} and \textnormal{`nonabelian'} will refer to the group structure of a quantum group, and \textnormal{`commutative'} and \textnormal{`noncommutative'} will refer to its topology.
\end{defn}
With this definition, the main theme of the paper is to simply replace $C^*$-algebras by the corresponding quantum spaces, and deduce results that can not be deduced, or even formulated, if we simply think of $C^*$-algberas as some abstract ``quantum spaces".  As examples of such results, we present the following: \newline
$\mathbf{(1)}$ An extension  of the spectral theorem to \emph{arbitrary} bounded operators on a  Hilbert space, i.e. an infinite dimensional analog of Jordan canonical form (Theorem \ref{spectral}). The Gelfand duality for the \emph{commutative} $C^*$-algebra generated by a normal operator $a$ leads to the spectral theorem \cite{dunford}. When $a$ is not normal, our noncommutative Gelfand-Naimark duality  yields: (i) Infinite dimensional Jordan canonical form, generalizing the spectal theorem, and (ii) A noncommutative functional calculus for $a$ (Theorem \ref{functional calculus}). \newline
$\mathbf{(2)}$ A general existence theorem for invariant subspaces (Theorem \ref{invariant}): It is well known that compact operators and normal operators on a several-dimensional complex Hilbert space have nontrivial invariant subspaces. These cover abitrary operators in finite-dimensional case, and for uncountably many dimensions, the result holds almost trivially.  For the case of \emph{separable} infinite dimensional Hilbert space, the result has been extended from normal operators to increasingly larger classes of operators \cite{kubrusly}. Our result (Theorem \ref{invariant}) covers any nonzero operator on any complex Hilbert space with several dimensions.\newline
$\mathbf{(3)}$ An extension of Pontryagin duality to nonabelian groups (Theorem \ref{quantum pontryagin}), which inevitably includes quantum groups as defined in Definition \ref{quantum space}. The classical Pontryagin duality asserts that the dual $\widehat{G}$ of an abelian locally compact group $G$ is a locally compact abelian group, and $G\cong \widehat{\widehat{G}}$.  There are several approaches to extending Pontryagin to possibly nonabelian locally compact groups \cite{enock,enock1,hopf1,hopf2}. One such approach uses Hopf-von Neumann algebras \cite{enock}; another, equivalently \cite{enock1,hopf1}, uses Hopf-$C^*$-Algebras.  These approaches embed the category of locally compact groups into certain categories of bialgebras, formulate the duality there, and characterize the bialgebras coming from groups and their dual bialgebras. Instead, we identify and emphasize a quantum group $\mathbf{\widehat{G}}$ (See Definition \ref{quantum space}) as the noncommutative  \textbf{dual of $\mathbf{G}$}. Then, the double dual $\widehat{\widehat{G}}$ is a locally compact group naturally homeomorphically isomorphic to $G$. Since quantum groups already appear in this theorem, we present an extension of Pontryagin duality to all quantum groups (Theorem \ref{quantum pontryagin}). This duality is entirely topological, and does not assume or employ any special Haar measure.  \newline
$\mathbf{(4)}$ The final section mentions two more applications: (i) Quick proof of Dauns-Hoffman theorems.  (ii) An extension of Stone's representation of Boolean algebras to orthomodular lattices. The second one applies the main idea of Theorem \ref{main} to an analogous problem in the field of Orthomodular lattices.

\section{\normalsize{NONCOMMUTATIVE\ \ GELFAND-NAIMARK\ \ DUALITY}}\label{mainsection}
Let $A$ be a $C^*$-algebra. A state $\alpha$ of $A$ is pure if and only if the corresponding Gelfand-Naimark-Segal (GNS)  representation $\pi_{\alpha}$ is irreducible \cite{kadison}. For $A$ without unit, $\overline{PS(A)}$, the weak$^*$-closure of the set of pure states, contains $0$. Then $P(A) := \overline{PS(A)}- \{0\}$ is locally compact Hausdorff, and is compact if and only if $A$ is unital.

\begin{defn}\label{equiv}
$\alpha , \beta \in P(A)$ are called \textbf{equivalent} if the corresponding GNS representations are equivalent. This is equivalent to saying that $\exists$ a unitary $u\in A$ such that for all $a\in A$, $\alpha(a) = \beta(uau^*).$
We denote this equivalence relation by $R(A)\subset P(A)\times P(A)$.
\end{defn}
\begin{proposition} 
A $C^*$-algebra $A$ is commutative if and only if the equivalence relation $R(A)$ is discrete, i.e. all its equivalence classes are singleton sets.
\end{proposition}
\begin{defn}\label{alpha_a}
For $\alpha\in P(A)$ define $\alpha_{a}$ by
\[\alpha_{a}(x) := \alpha(ax), \forall x\in A.\] 
\end{defn}
Observe that when $A$ is commutative, $\alpha_{a}(x) = \alpha(a)\,\alpha(x)$. Now, $\alpha_a$ can be uniquely extended to all bounded linear functionals by convexity, linearity and continuity.

\begin{defn}\label{a-hat}
For each $a\in{A},$ define a bounded linear operator \[\widehat{a}: C_0(P(A))\rightarrow C_0(P(A))\]
\[\,\,\,\text{by}\,\,\,\
\widehat{a}(\phi)(\alpha) := \phi(\alpha_a), \ \ \text{for each}\ \ \phi\in C_0(P(A)).\]
\end{defn}
When $A$ is commutative, $\widehat{a}$ is the multiplication operator, $\widehat{a}(\phi)(\alpha) = \alpha(a)\,\phi(\alpha) = \widehat{a}(\alpha)\,\phi(\alpha),$ where $\widehat{a}(\alpha) := \alpha (a).$
In the general case, $R(A) =  \cup_{a\in A} Supp(d\widehat{a}).$ Let $C_0(X)$ be the $C^*$-algebra of continuous complex functions vanishing at infinity on a locally compact (Hausdorff) space $X$. Then the double dual of $C_0(X)$ is a von Neumann algebra, and its maximal ideal space $Y$ carries a canonical class of measures. What follows is independant of a choice of measure $m'$ in this class. Let $m$ be the image of $m'$ under the canonical onto map $Y\rightarrow X$. Then $m$ gives an embedding $C_0(X)\hookrightarrow M(X): \phi\mapsto\int_{X}\phi\,dm,$ where $M(X)$ is the Banach space of complex valued regular Borel measures on $X$. 

\begin{defn}\label{operator-measure}
For a bounded linear operator $\widehat{a}: C_0(X)\rightarrow M(X)$, define a canonical regular complex valued Borel measure $d\widehat{a}$ on $X \times X$  by the identity  \[\int{(f\otimes g)\hspace{1pt}d\widehat{a}}\,\,=  \int{g\ d (\widehat{a}f\!).}\] 
\end{defn}

\begin{defn}\label{contmeasure}
A measure $\mu \in M(X\times X)$ will be called \textbf{continuous}if $\mu = d(\widehat{a})$, for a linear operator $\widehat{a}: C_0(X)\rightarrow M(X)$ the image of which is contained in $C_0(X)\subset M(X)$. Let $CM(X\times X)$ denote the algebra of such measures. \end{defn}
Now let $X = P(A)$ defined \textit{supra}, $CM(R(A))$ the algebra of continuous measures on $P(A)\times P(A)$ with support contained in $R(A)$, and $CO(R(A))$ the corresponding algebra of bounded operators $C_0(P(A))\rightarrow C_0(P(A))$. For $a\in A$, let $\widehat{a}$ be the operator defined in Definition \ref{a-hat}, and let $d\widehat{a}$ be the corresponding measure as in Definition \ref{operator-measure}. Now, our main result is:

\begin{theorem}[Noncommutative Gelfand-Naimark]\label{main}
Let $A$ be a $C^*$-algebra, and $X$ a locally compact Hausdorff space. Then,

\begin{enumerate}
\item

The assignments $a\mapsto \widehat{a}\mapsto d\widehat{a}$ give the following $C^*$-isomorphisms:
\[A\hspace{4pt}\cong\hspace{4pt}CO(R(A))\ \ \cong\ \ CM(R(A)).\]

\item 
If $R(X)$ is an equivalence relation on $X$, then we have natural isomorphisms
\[R(X)\ \  \cong\ \  R(CO(R(X)))\ \ \cong\ \ R(CM(R(X)).\]

\end{enumerate}
\end{theorem}
The proof uses a well-known noncommutative generalization of the Stone-Weierstrass theorem \cite{dixmier}.
For commmutative $A$,
$P(A)$ is the maximal ideal space of $A$, $R(A)$ is the diagonal of $P(A)\times P(A)$, and hence
$CO(R(A)) = CM(R(A)) =  C_0(P(A))$, so we recover the Gelfand-Naimark theorem (Theorem \ref{gelfand}).

\section{\normalsize{INFINITE \ JORDAN \ CANONICAL \ FORM: Extended Spectral Theorem}}\label{jordansection}

In this section we present an extension of Jordan canonical form to infinite dimensional Hilbert spaces. As such, it will also be an extension of the spectral theorem to non-normal operators. 

Let  $a$ be a bounded operator on a Hilbert space $H$. Let $A$ be the unital $C^*$-algebra generated by $\{1,a\}$,  and $R(A)$ be the equivalence relation defined by $A$ on $P(A)$.   When $a$ is normal, $A$ is commutative, and $R(A) = P(A) = \sigma(a),$ the spectrum of $a$, and the spectral theorem says that $a =  \int_{\sigma (a)} z \, dE ,$ where $z$ is the inclusion $\sigma(a)\!\hookrightarrow\!\mathbb{C},$ and $E$ the spectral measure corresponding to $a$ \cite{dunford}.

In the general case, where $a$ is not assumed normal, it turns out that the corresponding formula is equally simple (Theorem \ref{spectral}). As is often the case, the main effort goes into identifying the right concepts. We first  look at the situation heuristically, motivating the precise formulation that follows it.

Consider the map $\widehat{z} : P(A)\rightarrow \sigma(a)$ given by \,$\widehat{z}(\alpha) := \alpha(a)$, and  the canonical quotient map $q\!:\!P(A)\!\rightarrow\!Sp(A)$. Let $\pi$ be the partition of $P(A)$ generated  by those given by $\widehat{z}$ and $q.$  Denote the corresponding equivalence relation by $R(a)\!=\!Supp(d\widehat{a})\!\subset\! R(A)$, and let $r\!:\!R(a)\!\rightarrow\!Y$,  \,$\widehat{r}\!:\!P(A)\!\rightarrow\!Y$,  and $\tau\!:\!Y\!\rightarrow\!\sigma(a)$ be the natural maps:
\[
\xymatrix{
P(A)\ar@{->}[r]^{\widehat{z}}\ar@{->}[d]_{\widehat{r}}&\sigma(a)& R(a)\ar@{->}[r]^{z}\ar@{->}[d]_{r}  &\sigma(a)\\
Y\ar@{->}[ru]_{\tau}&& Y\ar@{->}[ru]_{\tau}}
\]
We think of\,$\tau\!:\!Y\!\rightarrow\!\sigma(a)$ as a uniformization of $\sigma (a)$: Corresponding to each $\lambda\in \sigma (a)$, there can be  several $y \in \tau^{-1}(\lambda) \subset Y,$  and for each such $y$,  there is an $R(a)$-block.  Thus, for each $\lambda\in \sigma (a)$ there are several $R(a)$-blocks. This is analogous to the finite dimensional case, where, for a fixed $\lambda\in\sigma(a)$, we may have several Jordan companion matrices $J_\lambda^k,$ of several different ranks $k,$ filling several disjoint diagonal square blocks of the Jordan canonical form.  We think of $R(a)$ as the scheme of blank blocks, to be filled with `Jordan matrices'. Before we make this rigorous by generalizing the notion of a spectral measure, we identify the algebra of sets on which it will be defined:

\begin{defn}[Lattice $\Omega$]
Let $R(X)$ be an equivalnce relation on a compact Hausdorff space $X.$ By \textbf{a sub-equivalence relation of $R(X)$}, we shall  mean  an equivalence relation $U$ on a subset of $X$, such that $U\subset R(X).$ We denote by $\Omega$ the set of all Borel sub-equivalence relations of $R(X)$. For $U, V \in \Omega ,$ define
 \begin{enumerate}
\item
$U \vee \, V\, :=$ The smallest sub-equivalence relation of $R(X)$ containing $U$ and $V$.
\item
$U \wedge \,V\, := \,\,U \cap \,V.$ 
\item
$1:= R(X),\, \, \, 0:= \emptyset $. 
\end{enumerate}
\end{defn}
Then $\Omega(\wedge, \vee , 1, 0)$ forms a lattice, which is not distributive in general. Consider the relational product $U \circ V := \{(x, y)\!\in\!R(X): \exists\,z\!\in\!X, (x, z)\!\in\!U, (z, y)\!\in\!V\}.$ Then all $U\in\Omega$ are idempotents, $U\circ U = U$. Also, $U\vee V = U \circ V$ if and only if $U\circ V = V \circ U$, in which case, we say that $U$ and $V$ \textit{commute}. It can be shown that $\Omega$ is distributive if and only if all $U, V \in\Omega$ commute.  Let $P(H)$ be the lattice of projections on a Hilbert space $H$.
Now we are ready to define the central notion of this section:

\begin{defn}\label{elementary measure}
An \textbf{Elementary measure on R(X) with respect to a Hilbert space $H$} is a function $E: \Omega\rightarrow P(H)$
which satisfies the following conditions:
\begin{enumerate}
\item
$E(\emptyset) = 0,$ $E(R(X)) = 1.$
\item
$E(U\wedge V) = E(U) \wedge E(V), \ \ \text{for all}\ \ U, V \in \Omega.$ 
\item
For a sequence $\{U_n\in\Omega , n = 1, 2, 3, \ldots\}$, such that $U_i \wedge U_j = 0$, if $i\neq j$,
\[E(\bigvee_{n = 1}^{\infty}U_n) =  \bigvee_{n = 1}^{\infty} E(U_n).\]
\end{enumerate}
\end{defn}
This definition differs that of  a spectral measure $E$ of a normal operator, in which case condition 3. reads $E(U\wedge V) = E(U)E(V).$ This identity implies that the image of $E$ is a Boolean algebra of projections. Indeed, when $R(X) =  diag(X\times X) \cong X,$ then $\Omega$ is the Borel algebra of $X,$ which is a Boolean algebra, its image under $E$ is a Boolean algebra of projections, and $E$ is simply a spectral measure.

Let $ a \in B(H)$, and $A$ be the $C^*$-subalgebra of $B(H)$ generated by $\{1,a\}$. Let $R(A)$, $R(a)$ and $\Omega$ be as defined above. Let $CM(R(A)) \cong A\subset B(H): \mu \mapsto \mu(a),$ be the inverse of the isomorphism given by Theorem \ref{main}.  Let $\mu_z := d\widehat{a}$  so that $\mu_z(a) = a$. Then we have the following:

\begin{theorem}[Jordan Canonical Form]\label{spectral}
With the notation established above, there exists on $\Omega$ a unique elementary measure $E$ with respect to the Hilbert space $H$  such that
\begin{enumerate}
\item If \ $U$ contains a nonempty open set, $E(U) \neq 0.$
\item$E(U)\,a\,E(U)= a\, E(U) \ \text{for all}\ U\in\Omega .$
\item \begin{equation} a  = \int_{R(A)}\!\mu_z  * dE. \label{jordan}\end{equation}
\end{enumerate}
Furthermore, for each $R(a)$-block $U\in \Omega$, $S(U) = \mu_z(U) - z(U)E(U)$ is either zero or a stable co-isometry, that is, $S(U)^*$ is an isometry, and $S(U)^n \rightarrow 0$ as $n \rightarrow \infty$, and hence $S(U)$ is a backward shift operator.
\end{theorem}
Thus, when $a$ is an (arbitrary) operator on finite dimensional $H$,  we recover the Jordan canonical form.
On the other hand, when $a$ is normal on arbitrary Hilbert space $H$, $A$ is commutative, $R(a) = R(A) = P(A) = Y = \sigma (a)$, so that $E$ and $\mu_z$ are the spectral measure and  the identity function respectively on $\sigma(a)$, and the theorem reduces to the Spectral Theorem. 
The formula \ref{jordan} for $a$ in Theorem \ref{spectral} is a special case of a formula for the the functional calculus $CM(R(A))\cong A$. Indeed, we have a larger functional calculus analogous to the functional calculus $L^\infty (\sigma (a), \nu) \cong  {a}''\subset B(H)$ of a normal operator:

\begin{theorem}[Functional calculus for a bounded operator]\label{functional calculus}
Let  ${a}''$ be the von Neumann algebra generated by $a,$ and let $LM(R(A))$ be the von Neuman algebra generated by $d\widehat{a}$.
Then the functional calculus $CM(R(A))\cong A$ can be extended to the functional calculus \[LM(R(A)) \cong {a}''\subset B(H): \mu\mapsto \mu(a),\] 
and is given by the following formula:
\begin{equation}\mu(a) = \int_{R(A)}\mu * dE,\label{calculus}\end{equation}
Also, this is the unique calculus on $LM(R(A))$ which satisfies the following:
\begin{enumerate}
\item $\mu_1(a) = 1,$ where $\mu_1$ is the measure corresponding to $1\in CO(R(A))$.
\item{$\mu_z(a) = a$}\label{1}
\item{$\mu\mapsto \mu(a)$ is an isometric monomorphism.}
\item{$\mu\mapsto \mu(a)$ extends the Riesz functional calculus.}
\end{enumerate}

\end{theorem}
An easy consequence of Theorem \ref{spectral} is,
\begin{theorem}[Invariant Subspace Theorem]\label{invariant}
Every bounded operator on a complex Hilbert space of dimension greater than one has a nontrivial invariant subspace.
\end{theorem}
The proof mimics that for normal operators \cite{kubrusly}.

\section{\normalsize{NONABELIAN\ \ PONTRYAGIN\ \ DUALITY}}\label{pontryagin}
Recall that the set $\widehat{G}$ of characters of a locally compact abelian group $G$ forms a locally compact abelian group and the celebrated Pontryagin duality theorem gives a natural isomorphism $G \cong \widehat{\widehat{G}}$.  We find that extending this theorem to nonabelian groups leads us to quantum groups as defined in Definition \ref{quantum space}: Given a locally compact group $G$ its dual is a quantum group $\widehat{G}$, which is a group if and only if $G$ is abelian. The classical dual of a possibly nonabelian $G,$ i.e. the set of equivalence classes of irreducible unitary representations of $G$,  is the quotient space corresponding to $\widehat{G}.$  In the abelian case, $\widehat{G}$ coincides with the classical dual.  This viewpoint inevitably leads to an extension of  the duality to \emph{quantum groups}. 

A quantum group $G$ as defined in Definition \ref{quantum space}  is a quotient $q: X\rightarrow Y$, the cartesian product $G\times G$ in this category is the fibred product $\pi:X\times_{q}X \rightarrow Y$, and the multiplication and the inversion are fibred maps $(m, b)$ and $(\iota , \beta)$ respectively:
\[
\xymatrix{
 X\times_{q}X\ar@{->}[d]\ar@{->}[r]\ar@{->}^{\pi}[dr]& X\ar@{->}^{q}[d]&X\times_{q}X\ar@{->}_{\pi}[d]\ar@{->}^{m}[r]& X \ar@{->}^{q}[d]& X\ar@{->}_{q}[d]\ar@{->}^{\iota}[r]&X\ar@{->}^{q}[d]\\
X\ar@{->}_{q}[r]&Y	&Y\ar@{->}_{b}[r]& Y& Y \ar@{->}_{\beta}[r] &Y}
\]
Then $X\times_{q}X \subset X\times X$ is the equivalence relation $R\subset X\times X$  given by the quotient map $q$. Instead of $m(x, y)$, we will write $xy$.
Note that $R = \bigcup_{y\in Y}( q^{-1}(y)\times q^{-1}(y)),$  and $m$ maps $q^{-1}(y)\times q^{-1}(y)\rightarrow q^{-1}(b(y))$ for each $y\in Y$. A quantum group $G$ is  a group if and only if $Y$ is a singleton set. On the other extreme, if $m$ is the second projection $R\rightarrow X$, then $G$ is simply the quantum space $X\rightarrow Y$, and when $X\times_{q}X  = diag(X\times X),$ then $G = X$ is merely a locally compact space.  Alternatively, we can view the multiplication as the partial map $R\times R\supset S\overset{m}{\longrightarrow} R$, given by $((x, y), (x', y'))\mapsto (xx', yy')$, where $S = \{((x,y), (x',y')): (x,y), (x',y'), (x,x'), (y,y') \in R\}.$ Let $CM(G)$ be the $C^*$-algebra of continuous Borel measures on $R$ as in Theorem \ref{main}, then $K^*(G) := CM(G)^{**}$  is a von Neumann algbera. Let $U, V$ be Borel subsets of $R$, and  $UV := \{m(u,v): (u,v)\in S, u\in U, v\in V\}.$ We define a co-multiplication $d: K^*(G)\rightarrow  K^*(G)\otimes K^*(G)$ by
\[d(\mu)(U\times V) = \mu(UV),\]
which makes $K^*(G)$ a  von Neumann-bi-algebra. We emphasize that $d$ may be a degenerate co-multiplication, and is nondegenerate if and only if $G$ is a group if and only if $X\times_q X = X\times X$. Also, the inversion map on $G$  gives an involution on $K^*(G)$. Thus, $K^*(G)$ satisfies all but one axioms of involutive Hopf-von Neumann algebras \cite{hopf2}.  We call such an algebra \textbf{a $\mathbf{K^*}$-algebra}. Then, the set $K^*(G)^*$ of weakly continuous functionals on $K^*(G)$ is naturally a Banach algebra. Let  $D^*(G)$ be the enveloping $C^*$-algebra of $K^*(G)^*$. Then, the \textbf{dual $\mathbf{K^*}$-algebra $\mathbf{\widehat{K^*(G)}}$} of $K^*(G)$ is defined  to be the enveloping von Neumann algebra of $D^*(G).$ Then $\widehat{\widehat{K^*(G)}} = K^*(G),$ and $D^*(G)$ is weakly dense in $\widehat{K^*(G)}.$ Now, using Theorem \ref{main},  we can construct from $D^*(G)$ a locally compact quantum space $\widehat{G}$, which has a multiplication  structure derived from the co-multiplication of $\widehat{K^*(G)}$.  This makes $\widehat{G}$ a quantum group which we call \textbf{the dual quantum group of $\mathbf{G}$}.  Following the same procedure, we construct a locally compact quantum group $\widehat{\widehat{G}}$ from $\widehat{K^*(\widehat{G})}$ and the main theorem of this section is:
\begin{theorem}[Pontryagin for quantum groups]
For a quantum group $G$,  $K^*(\widehat{G}) \cong  \widehat{K^*(G)}$,  $\widehat{K^*(\widehat{G})}$ $\cong  K^*(G),$ and $G \cong \widehat{\widehat{G}},$ and the following schematic diagram summarizes the situation:
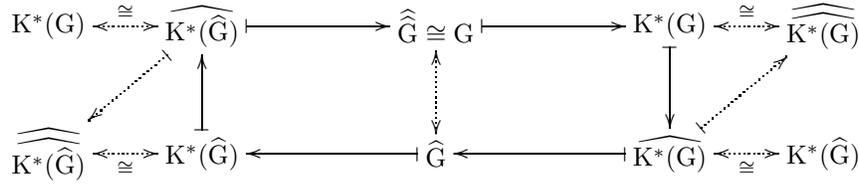
\[
\xymatrix{
K^*(G)\ar@{<.>}^{\cong}[r]&\widehat{K^*(\widehat{G})}\ar@{<-|}[d]\ar@{|->}[rr]&& {\widehat{\widehat{G}}}\cong \ar@{<.>}[d]G \ar@{|->}[rr]
&&  K^*(G)\ar@{|->}[d]\ar@{<.>}^{\cong}[r]& \widehat{\widehat{K^*(G)}}\ar@{<.|}[ld]\\
\widehat{\widehat{K^*(\widehat{G})}}\ar@{<.|}[ru]\ar@{<.>}_{\cong}[r]&K^*(\widehat{G})\ar@{<-|}[rr]&&\widehat{G} \ar@{<-|}[rr]&&\widehat{K^*(G)}\ar@{<.>}_{\cong}[r]&K^*(\widehat{G})
}
\]
\label{quantum pontryagin}
\end{theorem}
When $G$ is a locally compact group, $\widehat{G}$ is an abelian  \emph{quantum} group. Furthermore, if $G$ is an abelian locally compact group, then $\widehat{G}$ is an abelian locally compact group. Thus, the classical Pontryagin duality
is subsumed under the vertical arrows in the above diagram. On the other hand, when $G$ is just a quantum space, we have $K^*(G) \cong \widehat{K^*{G}} \cong  K^*(\widehat{G}) \cong  \widehat{K^*(\widehat{G})}$, with trivial co-multiplication, so that $G \cong \widehat{G}\cong \widehat{\widehat{G}}.$ That is, a quantum space is self-dual in this setting,
and the diagram reduces to only two arrows: $G\mapsto K^*(G)$ and $K^*(G)\mapsto G$. Thus, the horizontal arrows of the diagram subsume, via $D^*(G)$, the generalized Gelfand duality (Theorem \ref{main}).

Recall from Definition \ref{quantum space} that the terms `abelian' and `nonabelian'  refer to the group structure of a quantum group $G$, and `commutative' and `noncommutative' refer to the topology of $G$. Now,
let $G, H, K, N$ be quantum groups with the corresponding duals $\widehat{G}, \widehat{H}, \widehat{K}, \widehat{N}.$ Then the following table summarizes the various situations covered by Theorem  \ref{quantum pontryagin}: 

\begin{center}\begin{tabular}[h]{|c|c|c|}
\hline \boldmath{$\downarrow Topology$}\ \  \vline ${\ \ Group\rightarrow}$& Abelian & Nonabelian\\ \hline  Commutative & \small{$G, \ \ \overset{}{\widehat{G}}$}&  
\small{$K,\ \ \overset{}{ \widehat{H}}$}\\\hline Non-Commutative & \small{$H , \ \  \widehat{K}$}& 
\small{$N,  \ \ \overset{}{\widehat{N}}$}\\\hline
\end{tabular}
\end{center}
Thus, the dual $\widehat{G}$ of an abelian group $G$ is an abelian group;  for a nonabelian group $K$,  $\widehat{K}$ is an abelian noncommutative quantum group, etc.  We note that the box containing $G, \widehat{G}$ is the classical Pontryagin duality. The boxes containing $H, \widehat{K}$ and  $K, \widehat{H}$ include nonabelian groups and abelian noncommutative groups, and finally the box containing $N, \widehat{N}$ takes care of nonabelian noncommutative quantum groups.

\section{\normalsize{MISCELLANEOUS\ \ COMMENTS}}
\textbf{(1)} \ Theorem \ref{main} yields quick proofs of Dauns-Hoffmann theorems  \cite{dauns-hoffman}: $(i)$ Representations of a $C^*$-algbera as sections of a canonical ``sheaf" of (presumably simpler) $C^*$-algberas. $(ii)$ The center of a $C^*$-algbera
$A$ is isomorphic to $C_0(X)$, where $X$ is the spectrum of $A$.\newline
\newline
\textbf{(2)} \ The ideas of \S \ref{mainsection}  can be applied to Orthomodular Lattices (OML):
\begin{defn}\label{OML}
A set $L$ with operations $(\wedge , \vee , ', 0, 1)$  \textbf{is an orthomodular lattice} if  $\forall s, t \in L, \ (i) \,L(\wedge , \vee)$ is a lattice,  $(ii)$\, $(s')' = s ,$\ 
$(iii)$\, $s \leqslant t  \Longrightarrow   t' \leqslant s' ,$ \ $(iv)$ \,$ s \vee s' = 1,  s \wedge s'  =  0,$\  $(v)$\, $ s 
\leqslant t \Longrightarrow s \vee (s' \wedge t) = t$.
Condition (v) is called the \textbf{orthomodularity} condition. 
\end{defn}
Lattices of projections in a $C^*$-algebras are the prime examples of OML.
Note that orthomodularity is a weakening of distributivity, so that distributive OML are simply Boolean algebras. If for $a, b \in L$ we set $a\ \Dot{\wedge}\  b := (a\vee b')\wedge b,$ then $L$ is a Boolean algebra if and only if\ $\forall \, a, b \in L, a \ \Dot{\wedge}\  b = b\  \Dot{\wedge}\  a$ \cite{beran}, in which case, $a \  \Dot{\wedge}\  b = a \wedge b.$ In this sense, an OML is a noncommutative generalization of Boolean algebra. Also, if a $C^*$-algebra $A$ is generated by its lattice $L_A$ of projections (for example when $A$ is a von Neumann algebra), then $A$ is commutative if and only if the OML $L_A$ is commutative, i.e. a Boolean algebra. Now, elements of a commutative OML, i.e. a Boolean algebra, are represented by clopen sets of a totally disconnected compact space---its maximal ideal space (Stone's Theorem \cite{stone}).  As in the case of $C^*$-algebras, the geometric object corresponding to a (possibly noncommutative)\, OML is an equivalence relation on (or a quotient of) a totally disconnected compact space naturally associated with the lattice. Furthermore, an OML is Boolean if and only if this equivalence relation is discrete. In this case, one recovers Stone's theorem. The general case yields Graves-Selesnick representation \cite{selesnick} of an OML as sections of sheaves of (presumably simpler) OML's,  an OML analog of Dauns-Hoffman theorems.

\end{document}